\documentclass[11pt]{amsart}
\usepackage{amsfonts}
\usepackage{amsmath}
\usepackage{amssymb}
\usepackage{graphicx}%
\newtheorem{theorem}{Theorem}[section]
\theoremstyle{plain}

\newtheorem{lemma}{Lemma}[section]

\newtheorem{proposition}{Proposition}[section]

\numberwithin{equation}{section}

\usepackage{amsmath}
\usepackage{verbatim}

\usepackage{enumerate}
\usepackage{amsmath}
\usepackage{amsfonts}

\newcommand{\R}{{\mathbb{R}}}

\newcommand{\ds}{\displaystyle}

\newcommand{\Rnn}{{\R^n \times \R^n}}

\newcommand{\D}{{\nabla}}

\newcommand{\C}{{\mathbb{C}}}
\newcommand{\I}{{\sqrt{-1}}}
\newcommand{\N}{{\Sigma}}
\newcommand{\remove}{{\setminus}}

\newcommand{\thetat}{{t}}

\begin{document}
\title[Calibrations associated to Monge-Amp\`ere equations]{Calibrations associated to Monge-Amp\`ere Equations}
\author{Micah WARREN}
\address{Department of Mathematics, Box 354350\\
 University of Washington\\
Seattle, WA 98195}
\email{mwarren@math.washington.edu}
\thanks{}
\date{July 9, 2006}
\begin{abstract} We show the volume maximizing property of the special Lagrangian submanifolds of a pseudo-Euclidean space. These special Lagrangian submanifolds arise locally as gradient graphs of solutions to Monge-Amp\`ere Equations.  
\end{abstract} 
\maketitle

\section[Intro]{Introduction}

In this note we show

\begin{theorem}
Let $ \Gamma = (x, \D u(x) ) \subset \Rnn $ be the gradient graph of a convex function $u\in C^2 (\Omega),$ for $\Omega$ a bounded, simply connected region with $C^1$ connected boundary. If $u$ satisfies the Monge-Amp\`ere equation \begin{equation}\det (D^2u) = c, \end{equation} or equivalently,  $$ \ds\sum \ln \lambda_i = \ln c, $$ where $\lambda_i$ are the eigenvalues of $D^2u,$  then $\Gamma$ is volume maximizing in the pseudo-Euclidean space  $(\R^n_x \times \R^n_y,dx dy). $  Precisely, if $\Sigma$ is any $C^1$, space-like, oriented $n$-surface, with $\partial \Sigma  = \partial \Gamma = \Gamma \cap \{x \in \partial \Omega\},$  homologous to $\Gamma $ in $(\R^n \times \R^n, dx dy ) ,$ then $$Vol(\Sigma) \leq Vol(\Gamma) $$ with equality only if $\Sigma = \Gamma$.


\end{theorem}

A Lagrangian submanifold of $\Rnn$ is one that can be described locally as a gradient graph, $(x,\nabla u(x)).$  The metric $ dxdy$ can be expressed explicitly as the indefinite form $$ dxdy =   \frac{1}{2} \sum_i (dx_i \otimes dy_i + dy_i \otimes dx_i).$$ Hitchin [H, \textsection 5] introduced a definition of ``special Lagrangian"  for Lagrangian submanifolds  of  $(\R^n \times \R^n, dxdy) ,$ and demonstrated that a gradient graph $(x, \D u(x))$ is special precisely when the potential $u(x)$ satisfies (1.1).  Jost and Xin [JX, \textsection 4] then showed that such a submanifold has mean curvature $ H \equiv 0.$  We now show that if  $u(x)$ is a solution to (1.1), then the gradient graph $(x, \D u(x))$ is a calibrated submanifold of $(\R^n \times \R^n, dxdy),$ and consequently volume maximizing.  In fact, any submanifold which is locally described by gradient graphs of functions satisfying (1.1) is calibrated, and therefore volume maximizing.   Theorem 1.1 can be stated to allow a slightly larger class of $n$-surfaces, as we shall see in section 3.

\bigskip

The study of calibrated Lagrangian submanifolds of $\Rnn$ began with the work of Harvey and Lawson [HL, \textsection 3], who studied Lagrangian submanifolds of  $ \C^{n} \cong \Rnn$ with the Euclidean metric $\delta_0$, and showed that a Lagrangian submanifold is calibrated, and therefore volume minimizing, if and only if the potential $u(x)$ is a solution to the special Lagrangian equation 

\begin{equation}  \sum_i \arctan \lambda_i = c \,. \end{equation}

\bigskip

By taking linear combinations of the metrics $\delta_0$ and $g_0 = 2 dxdy$, we obtain a family of metrics on $\Rnn$

  $$g_\thetat  = \cos\thetat \,g_0 + \sin\thetat \,\delta_0 .$$
\bigskip

We find that the extremal Lagrangian surfaces in $(\Rnn, g_{\thetat})$ arise as solutions to a family of special Lagrangian equations

\begin{equation} For \,\, \thetat = 0 \quad \sum_i \ln \lambda_i = c \end{equation}

\begin{equation} For \,\, \thetat \in (0,\frac{\pi}{4})  \quad \sum_i \ln(\frac{\lambda_i+ a -b }{\lambda_i+ a+b} )= c\,\end{equation} 

\begin{equation} For \,\, \thetat = \frac{\pi}{4} \quad \sum_i \frac{1}{1+ \lambda_i} = c \end{equation}

\begin{equation}For \,\,  \thetat \in (\frac{\pi}{4}, \frac{\pi}{2})  \quad \sum_i \arctan ( \frac {\lambda_i+ a -b } {\lambda_i+ a + b}) = c \,\end{equation}
 
\begin{equation}For \,\,  \thetat =  \frac{\pi}{2} \quad \sum_i \arctan  \lambda_i  = c \,\end{equation}where  $a= \cot \thetat$ and $b = \sqrt{|\cot^2 \thetat - 1| } .$  Further, we have the following extremal volume property of special Lagrangian graphs in $(\Rnn, g_\thetat)$

\begin{theorem}

i) Suppose $u \in C^2(\Omega),\,\, \Omega \subset \R^n,$ is a solution to (1.4).    If the metric  $g_\thetat$ restricts to a positive definite metric on $ \Gamma = (x, \D u(x) ), $ then $\Gamma$ is volume maximizing among homologous, $C^1$, space-like $n$-surfaces in $ (\Rnn, g_\thetat) ,$ as in Theorem 1.1. 

ii) Suppose $u \in C^2(\Omega),\,\, \Omega \subset \R^n,$ is a solution to (1.5). Then the volume of $\Gamma = (x, \D u(x) ) $ is equal to the volume of any homologous, $C^1,$ space-like $n$-surface in $ (\Rnn, g_{\frac{\pi}{4}})$. 

iii)  Suppose $u \in C^2(\Omega),\,\, \Omega \subset \R^n,$ is a solution to (1.6). Then  $\Gamma = (x, \D u(x) ) $ is absolutely volume minimizing among all homologous $n$-surfaces in $ (\Rnn, g_\thetat)$.

\end{theorem}

Using a change of variables, we may restate the following Bernstein-type results of Calabi-Pogorelov, Flanders, and Yuan

\begin{theorem} 

i)\normalfont (Calabi [C], Pogorelov [P]) \it Suppose  $u \in C^2(\R^n)$, $D^2 u \geq -\cot \thetat $, and $\Gamma = (x, \D u(x))$ is a maximal space-like surface in $(\Rnn, g_\thetat) $.  Then $u(x) $ is a quadratic polynomial.   

ii)\normalfont (Flanders [F]) \it Suppose $u \in C^2(\R^n) $ is a convex solution to (1.5). Then $u(x)$ is a quadratic polynomial.  

iii)\normalfont (Yuan [Y1],[Y2]) \it 
  There is a value $C_\thetat$ such that if $u \in C^2(\R^n) $ is a solution to (1.6), with either

a) $D^2 u  \geq - \cot \thetat $ or 

b) $ c > C_\thetat + \frac{(n-2)\pi}{2}$,  

\noindent then $u(x) $ is a quadratic polynomial.  

\end{theorem}

In section 2 we give a brief background of the related Lagrangian and calibrated geometry. In section 3 we show that calibrations associated to Monge-Amp\`ere equations detect maximal surfaces, proving Theorem 1.1.  We will derive the special Lagrangian equations using mean curvature in section 4, and in section 5, using a Lewy rotation, we prove Theorem 1.2 and the related Bernstein results.  In section 6,  we provide a simple example of a maximal Lagrangian surface in $(\R^2\times\R^2, g_0)$ which cannot be globally described by a solution to (1.1),  but which is described by a solution to (1.4). 

\bf Acknowledgment. \normalfont The author is happy to thank to his thesis advisor, Yu Yuan, for suggesting this problem, and for the many fruitful discussions about this and other problems.

\section[]{Background}

Let $\xi$ be an oriented real $n$-plane in $\Rnn \cong \C^n$.  If the projection $\Rnn \rightarrow \R^n_x$ is full rank on $\xi,$ we may write $\xi$ as the span of the vectors

$$ \xi_i = \partial_{x_i} + w_i^{\,\,j}  \partial_{y_j}\quad i=1,...,n.$$

In this form, it is easy to see that $J \xi \perp  \xi$ if and only if $ w_j^{\,\,i} = w_i^{\,\,j},$ where $J$ is the automorphism giving $\C^n$ its standard complex structure.  We say an $n$-plane is   \it Lagrangian \normalfont if it has the property $J \xi = \xi ^\perp, $ and call the set of such planes $L.$

Let $\Gamma = (x, F(x))$ be the graph of an $\R^n$-valued function in $\Rnn$.   We say the submanifold $\Gamma$ is Lagrangian if, at each point $p$, $T_p\Gamma$ is Lagrangian, or equivalently, \begin{equation}\frac{\partial F^i}{\partial{x_j}} =  \frac{\partial F^j}{\partial{x_i}} \end{equation} which is the case precisely when $F  = \nabla u $, locally, for some $u.$ 

For $\Gamma$ the graph of $F$ over a bounded region $\Omega \subset \R^n$, the volume of $\Gamma$ is given by the functional  

$$ Vol(\Gamma) = \int_\Omega \sqrt{ \det( I + (DF)(DF)^T) }dx  $$ when the metric on $\Rnn $ is the Euclidean metric $\delta_0.$ When the ambient metric is the pseudo-Euclidean metric $g_0,$ the volume functional is given by

$$ Vol(\Gamma) = \int_\Omega \sqrt{\det( DF + (DF)^T)}dx \,.$$ 

When $F = \D u $, then the Euler-Lagrange equations for the above functionals at $F$ become (1.2) and (1.1), respectively.

 We are interested, however,  in showing that the solutions to these Euler-Lagrange equations provide not only critical points,  but in fact global extrema for the volume functional, over all homologous $n$-surfaces.  Furthermore, a given Lagrangian surface may arise as a graph only locally, so we desire a way to compare non-graphical Lagrangian submanifolds to their homologous competitors, in a way which only uses local information about the surface.  These concerns motivate the use of a calibration.

We recall from [HL,  \textsection 1] the notion of a calibration.  Given a Riemannian manifold $M$, suppose $\varphi$ is a closed exterior $k$-form on $M$ with the property that on all oriented tangent $k$-planes $\xi$,

$$ \varphi (\xi) \leq Vol(\xi)\,.  $$ We call such an $M$ a \it calibrated manifold\normalfont. If $\N \subset M $ is a compact, oriented, $k$-dimensional submanifold with the property that 

\begin{equation} \varphi|_\N(\xi) = Vol|_\N(\xi) \, ,\end{equation}  for all oriented (with respect to the orientation on $\N$)  tangent $k$-planes $\xi$, then $\N$ is homologically volume minimizing, and we call $\N$ a \it calibrated submanifold \normalfont.   Recall (cf [L, p. 431]) that $\N$ is homologous to $\N'$ if there exist smooth triangulations $c$ and $c'$ of $\N$ and $\N'$ such that $c-c'$ is a boundary.  If $ \N$ is homologous to $ \N' $,  an application of Stokes's Theorem gives 

$$ Vol(\N)\, =\, \int_\N \varphi\, = \,\int_{\N'} \varphi\, \leq \, Vol(\N')\,.$$

We may apply this idea in the pseudo-Riemannian setting, where we are looking for volume maximizing submanifolds, provided we pay careful attention to orientation as it arises.  Given an oriented $k$-plane $\xi$ and a $k$-form $\varphi,$ we will say that $\xi$ is \it oriented with respect to $\varphi$ \normalfont whenever $\varphi(\xi) > 0.$  Let $M$ be a pseudo-Riemannian manifold with index $n$$-$$k$, and suppose $\varphi$ is a closed exterior $k$-form on $M$ with the property that on all oriented space-like $k$-planes, which are oriented with respect to $\varphi,$  
$$ \varphi (\xi) \geq Vol(\xi). $$ Let $\N$ be any compact, oriented, $k$-dimensional, space-like submanifold with the property that $\varphi|_\N(\xi) = Vol|_\N(\xi)$, for all oriented (with respect to the orientation on $\N$) tangent $k$-planes $\xi$.  If $\N$ is homologous to $\N'$, then $\N'$ has a prescribed orientation, and it is with respect to this orientation that we apply Stokes's Theorem and obtain $$ \int_\N \varphi\, = \,\int_{\N'}\varphi .$$ It is possible, however, for the orientation of $\N'$ to produce tangent planes which are not oriented with respect to $\varphi.$ If this disagreement occurs, the calibrating inequality is reversed, and we are unable to make the volume comparison. This can happen if $\Sigma'$ is disconnected, or is singular along a significant subset.  The conditions in Theorem 1.1 are sufficient to preclude any such pathology, as we will see in section 3.

\bigskip

For the Euclidean case, Harvey and Lawson [HL, \textsection 3] define the $S^1$-family of forms $$\alpha_\theta = Re(e^{-\I\theta}dz),$$ where $ dz = dz_1\wedge ...\wedge dz_n,$ and show that these are calibrations on $\C^n.$    We briefly recall the idea.

For any real $n$-plane $\xi $ in $\Rnn = \C^n,$ we may choose an oriented orthonormal basis $\xi_1,...,\xi_n,$ with $\xi_i \in \C^n$.  Define a complex-linear map $ A: \C^n \rightarrow \C^n$ by $$A(e_{x_i}) = \xi_i ,$$ $$ A(e_{y_i}) = J \xi_i ,$$ which is represented by a complex-valued matrix $A \in M(n, \C),$ namely the complex $n \times n$ matrix with columns  $\xi_1,..., \xi_n.$  This also defines a real-linear map $ \hat A: \Rnn \rightarrow \Rnn,$ which is represented by a real-valued matrix $\hat A \in M(2n, \R) .$ 
Now  $$|\xi  \wedge J\xi | = \mbox{det}_\R \hat A = |\mbox{det}_\C A|^2 = |dz(\xi)|^2 = \alpha^2(\xi)+\beta^2(\xi),$$ where $\alpha(\xi) = Re(dz(\xi)),$ and $\beta(\xi) = Im(dz(\xi)). $ By Hadamard's Inequality,  $|\xi|^2 \geq |\xi\wedge J\xi|$, with equality if and only if $\xi$ is Lagrangian.   So now we have $$ |\xi|^2 \geq |\xi\wedge J\xi| = \alpha^2(\xi)+\beta^2(\xi) \geq \alpha^2(\xi),$$ with equality if and only if $\xi$ is Lagrangian and $\beta(\xi) = 0 .$

If $\xi$ is Lagrangian, then $\xi_1, ...,\xi_n$ is an orthonormal basis for $\C^n$, so $A \in U(n).$  With the action of $SO(n)$ on $U(n)$ given by $n \times n$ matrix multiplication, choosing an orthonormal basis associates to each $\xi \in L$ an $ A \in U(n)$  which is unique up to a factor of $SO(n)$. With this association, we have a transitive action of $U(n)$ on $L.$ The isotropy group of $\xi_0 =  \R^n_x$ is $SO(n)$, so $L$ is a homogeneous space and the complex determinant descends to a map
$$\Theta: L \approx U(n) / SO(n) \overset{\det_\C}\rightarrow S^1.$$  For any $\xi \in L$, $dz(\xi) = e^{\I\Theta(\xi)} $. Taking 
$$\alpha_{\theta} = Re(e^{-\I\theta}dz),$$ we see that $\alpha_{\theta} (\xi) = 1$ if and only if $\Theta(\xi) =  \theta,$ thus the $\alpha_\theta's$ are calibrations for $\C^n$.

The special Lagrangian equations may then be deduced from the condition $$ \arg (\det(I+\I D^2 u)) = \arg((1+ \I\lambda_1) ... (1+\I\lambda_n)) = c $$ for some fixed $c$. This equation is satisfied precisely when
$\Theta $ is constant along the graph of $\D u$.

McLean [M, Theorem 3.6, 3.10] showed that given $\N \subset M $ a compact special Lagrangian submanifold, the moduli space $X$ of special Lagrangian manifolds near $\N$ in $M$ is itself a manifold which carries a natural Riemannian metric.  Hitchin [H, Proposition 2], showed that this metric can be obtained locally by embedding the moduli space $X$ into $(V\oplus V^*, g_0),$  where $V = H^1(\N,\R)$.  Hitchin then showed [H, Proposition 3] that $X$ is a special Lagrangian submanifold of $(V\oplus V^*, g_0)$, and that special Lagrangian submanifolds arise as solutions to the Monge-Amp\`ere equation, $\det(D^2 u ) = c.$

\section[]{Calibrations for Pseudo-Euclidean Space }

\bigskip

In the pseudo-Euclidean setting, a Lagrangian submanifold $\Sigma$ is called \it special \normalfont [H, p. 510] if a linear combination of the volume forms $dx_1\wedge ...\wedge dx_n$ and $dy_1\wedge ...\wedge dy_n$ vanish along $\Sigma$.

\begin{proposition}  For $ c > 0$,  $$ \Phi_c = \frac{1}{2} [ c \,dx_1\wedge...\wedge dx_n + \frac{1}{c}dy_1\wedge ...\wedge dy_n ] $$  is a calibration for $ (\Rnn, dxdy)$. Suppose $\xi$ is an oriented space-like $n$-plane in $\Rnn$, with $\Phi_c(\xi) > 0. $ Then  $$\Phi_c(\xi) \geq  Vol(\xi)\,$$ with equality if and only $\xi$ is special Lagrangian, that is, if $$ dy_1\wedge ...\wedge dy_n(\xi) = c^2 dx_1\wedge...\wedge dx_n(\xi) .  $$ 
\end{proposition}

In the development of the Euclidean calibrations, Harvey and Lawson used Hadamard's Inequality to compare Lagrangian planes to non-Lagrangian planes. In order to prove Proposition 3.1, we need a result which serves this purpose in the pseudo-Euclidean case. We recall the following lemma from linear algebra.

\begin{lemma} Suppose $Q \in GL(n,\R)$ satisfies $ Q_{ij}x^ix^j >\, 0, \, $ for all $ 0 \neq  x \in \R^n.$  Then $$\det(Q) \geq \det(\frac{Q+Q^T}{2})$$  with equality if and only if $Q = Q^T.$

\end{lemma}

\it Proof of Lemma 3.1. \normalfont  The matrix $Q$ can be written $Q = S+A$, with $S$ symmetric, and $A$ antisymmetric.  Choose a basis so that $S$ is a diagonal matrix. With respect to this basis, 

$$ Q = 
\begin{pmatrix} \lambda_1 & a_{12} &...  & a_{1n}   \\ -a_{12}  & \lambda_2   & ...  & a_{2n} \\  ...  &... & ...&... \\ -a_{1n} & -a_{2n} & ...   & \lambda_n  \\ \end{pmatrix}. $$

Expand the determinant of Q, and group the terms according to the number of $\lambda_i$'s each terms contains.  For all $k \leq n $ define 
$$ P_k  =  \mbox{sum of all terms containing exactly} \,\, k\,\,   \lambda_i's\,. $$ We see that $P_n$ consists of one term, namely $\lambda_1 ... \lambda_n = \sigma_n(S), $  and that there are no terms with $(n-1)$  $\lambda_i$'s, so $P_{n-1} = 0. $  For $P_{n-k},\, k \geq 2$, we fix $i_1,  ..., i_{n-k}$, and look at the terms containing $\lambda_{i_1}\, ... \,\lambda_{i_{n-k}}. $ These occur as the determinant of a matrix, which after an orthogonal change of basis, looks like $$  
\begin{pmatrix} \lambda_{i_1} &  &   &  &  & \\  & ...   &  &   &  & \\   &     & \lambda_{i_{n-k}}& & & \\  &  &   &  0 &a_{j_1j_2} & ... \\  &  &  & -a_{j_1j_2}&0  & \\  & & & ... &   &\end{pmatrix} $$ with $j_1,..., j_{k} \notin \{ i_1,..., i_{n-k}\}.$ This determinant is the product of the determinants of a positive diagonal matrix and an antisymmetric matrix.  It follows that  $P_{n-k} \geq 0 .$ We also see that $$P_{n-2} = \sum_{i<j} a_{ij}^2\,\,\lambda_1 ...\hat\lambda_i \hat\lambda_j... \lambda_n $$ which is strictly positive unless $a_{ij} = 0,\, $ for all $ \,i,j.$   We conclude $$ \det(Q) = P_n + P_{n-2} + ... +  P_0  \geq P_n = \det(S)=\det(\frac{Q+Q^T}{2}) $$  with equality if and only if $Q=Q^T$. $\qed$

\bigskip

\it Proof of Proposition 3.1. \normalfont If $\xi$ is a space-like tangent plane, the projection onto $\R^n_x$ is full rank, and we can take a basis for $\xi$ of the form
$$ \xi_i = \partial_{x_i} + w_i^{\,\,j} \partial_{y_j}, \quad i=1,...,n. $$ Let $Q$ be the matrix given by $ Q_{ij} = w_i^{\,\,j} $.  If $g$ is the induced metric, then with respect to this frame, the tensor $g_{ij} = g(\xi_i, \xi_j)$  becomes $(Q+Q^T)_{ij}/2$, so 
$$ Vol(\xi_1\wedge ... \wedge\xi_n) =  \sqrt{ \det(\frac{Q+Q^T}{2})}  .$$Now  $$ \Phi_c(\xi_1\wedge ...\wedge \xi_n)  = {\frac{1}{2}}[c + \frac{1}{c} \det(Q)] \geq    \sqrt{\det{Q}}$$ with equality if and only if $\det(Q) = c^2$, so  
$$ \Phi_c(\xi_1\wedge ...\wedge \xi_n)  \geq  \sqrt{ \det(Q)}.$$ From Lemma 3.1 $$ \sqrt{\det(Q)} \geq \sqrt{ \det(\frac{Q+Q^T}{2})} =Vol(\xi_1\wedge ...\wedge \xi_n) $$with equality only if Q is symmetric, that is, if $\xi$ is Lagrangian.  Hence
$$ \Phi_c(\xi_1\wedge ... \wedge\xi_n)\geq Vol(\xi_1\wedge ... \wedge\xi_n) $$with equality if and only if $\det(Q)= c^2 $ and $Q = Q^T.\quad \qed$

\bigskip
We prove a more general result than Theorem 1.1. 
\bigskip

\begin{theorem}  Suppose $\Omega \subset \R^n $ is a bounded region, and $ u \in C^2(\Omega)\,\, $  is a convex solution to the Monge-Amp\`ere equation $$\det (D^2u) = c^2.$$  Then in the pseudo-Euclidean space  $(\R^n_x \times \R^n_y, dxdy),$  the gradient graph of $u(x) $,  $ \Gamma = (x, \D u(x) ),$ is volume maximizing  in the following sense:

Let $\Sigma$ be an oriented $n$-surface which is homologous to $\Gamma $ in $\R^n \times \R^n, $ with $\partial \Sigma  = \partial \Gamma = \Gamma \cap \{x \in \partial \Omega\}.$ If any of the following hold \begin{itemize}\item $\partial \Omega$ is connected and $\Sigma$ is $C^1$ and space-like, 
\item $\Sigma$ is connected, $C^1$ and space-like, or \item $\Sigma$ is $C^1$ and space-like except on a set $\Sigma_0$, which has $n$-dimensional Lebesgue measure zero, such that $ \Sigma\remove \Sigma_0 $ is connected, 
\end{itemize}
then $Vol(\Gamma) \geq Vol(\Sigma) . $ If $ \Omega$ is simply connected with $C^1$ connected boundary, and $\Sigma$ is $C^1,$ then equality holds only if $\Sigma = \Gamma.$

\end{theorem}

\it Proof of Theorem 3.1, Theorem 1.1.  \normalfont Suppose $u \in C^2(\Omega)$ with $\det(D^2 u ) = c^2 $ on $\Omega.$ Then $$ \Phi_c(\partial_1\wedge ...\wedge \partial_n) = {\frac{1}{2}}[c + \frac{1}{c} c^2] =  c ,$$

$$ Vol(\partial_1\wedge ...\wedge \partial_n) = \sqrt{\det(D^2u)} = c.$$ So $$Vol|_\Gamma =  \Phi_c |_\Gamma,  $$ and $\Gamma$ is a calibrated submanifold. For any homologous $n$-surface $\Sigma$, we know that (cf [L, p. 431]) $$ \int_\N \Phi_c = \int_\Gamma \Phi_c = Vol(\Gamma),$$ so our task is to show that each of the listed conditions imply that the oriented planes for $\Sigma$ are oriented with respect to $\Phi_c$, giving $$ \Phi_c|_\N \geq Vol|_\N.$$ 

We begin by showing that the first condition implies the second condition. Regarded as a linear map, $g_0 \in GL(2n, \R)$ has two eigenvalues, $+1$ and $-1$, so we can decompose $\R^{2n}$ into the eigenspaces corresponding to these two eigenvalues. Let $P^+$ be the the projection of $\Rnn$ onto the eigenspace corresponding to the eigenvalue $+1$.    If $\N$ is $C^1$ and space-like, the projection $P^+|_\Sigma \rightarrow \R^n_+$ must be full rank, so must be an open map.  It follows that each component of $\Sigma$ must have non-empty boundary. Since $\partial \Sigma = \partial \Omega$, $\partial \Sigma $ must be connected. Each component of $\Sigma$ intersects $\partial \Sigma $, so $\Sigma$ consists of a single component.  

It is clear that the second condition implies the third.  So now assume $\N$ is $C^1$ and space-like except on a set $\Sigma_0$, which has $n$-dimensional Lebesgue measure zero, and $ \Sigma \remove \Sigma_0 $ is connected. The integral $ \int_\N \Phi_c $ is positive, so for the induced orientation on $\N$, $\Phi_c(\xi) $ is positive for some oriented tangent plane $\xi_s = T_s\N,$ at some point $s \in \N$.   Since $Vol$ does not vanish on $\N \remove \N_0$, it follows from the $C^1$ assumption that $\Phi_c(\xi_s)  \geq Vol(\xi_s) > 0  $ for all $s \in \N \remove \N_0$. Hence $$Vol(\N) =  \int_\N dVol \leq \int_\N \Phi_c =  Vol(\Gamma). \qed$$ 

\it Uniqueness. \normalfont It is clear by  Proposition 3.1 that equality will only occur if $\Sigma$ is a special Lagrangian surface, locally described by gradient graphs of $\det D^2 u = c^2.$  In order to use the comparison principle, we first show that $\Sigma$ is globally described by a single graph over $\Omega.$

Let $P^x$ be the the projection of $\Rnn$ onto $\R^n_x,$ and let $\Omega_0 = P^x(\Sigma).$   
We observe that, due to the space-like condition, the projection $P^x$ is open on the interior of $\N, $ so $\partial \Omega_0 \subset P^x(\partial \Sigma) = \partial \Omega.$ 

Let $p = (x_0, \nabla u (x_0)) \in \partial \Omega = \partial \Sigma,$ for $x_0$ an extreme point of $\Omega.$  From the open condition on the map $P^x$ on $\Sigma$, it follows that $x_0$ is also an extreme point for $\Omega_0,$ and that the inward pointing normals for $\partial \Omega$ and $\partial\Omega_0$ must agree at $x_0.$  The regions $\Omega_0$ and $\Omega$ then must intersect on a non-trivial open set near $x_0.$  Using the openness of $P^x,$ together with the boundary condition $\partial \Omega_0 \subset \partial \Omega,$  it is then easy to check that $\Omega_0 \cap \Omega$ is relatively closed and relatively open as a subset of $\Omega$, hence $\Omega \subset \Omega_0.$  
 
Take a cover of $\Sigma^\circ$ by open sets $\Sigma_i, $ where $\Sigma_i = \{ (x, \nabla v_i) | x \in U_i \} , $ for $v_i$ solutions to $\det D^2 v_i = c^2$ on $U_i.$  Let $\hat U_1 = U_1$, then recursively define $\hat U_i = U_i \remove \hat U_{i-1}$, and define $\hat \Sigma_i = \{ (x, \nabla v_i) | x \in \hat U_i \}  .$ 
The disjoint union $\bigcup_i \hat U_i$ is then an open set that contains $\Omega,$ and the disjoint union $\bigcup_i \hat \Sigma_i $ is a subset of $\Sigma.$

Suppose that $\bigcup_i \hat \Sigma_i \neq \Sigma^\circ.$  Then there exists an open subset $U_{jk} \subset U_j \cap U_k,$ with $j < k ,$ such that $\Sigma_j$ and $\Sigma_k$ are both graphs over $U_{jk},$ but are disjoint.   Then $\Sigma$ must contain the disjoint union $$ \Sigma' = \bigcup_i \{ (x, \nabla v_i) | x \in \hat U_i \} \bigcup \{(x, \nabla v_k) | x \in  U_{jk} \}.$$  However, integrating over $\Sigma$ gives $$ \int_\Sigma dVol \geq \int_{\Sigma'} dVol \geq \int_\Omega c \, dx + \int_{U_{jk}} c\, dx > c|\Omega| = Vol(\Gamma) $$ contradicting the inequality in the  conclusion of Theorem 3.1. We conclude that $\Sigma^\circ$ = $\bigcup_i \hat \Sigma_i$ and $P(\Sigma^\circ) = \Omega.$ The gradients $\nabla v_j$ and $\nabla v_k$ must agree on all overlaps $U_j \cap U_k,$ so we can extend any of the $\R^n$-valued functions $\nabla v_j$ to an $\R^n$-valued function $F$ on all of $\Omega.$  We have assumed that $\Omega$ is simply connected, so $F = \nabla v,$ for some $v$ satisfying $\det D^2 v = c^2, $ on all of $\Omega.$  To compare $u$ with $v$, we note that $\nabla u = \nabla v $ is already fixed around the boundary of $\Omega,$  which is connected, so we may integrate $u$ and $v$ around the boundary and conclude that $ u$ and $v$ differ by a constant. Applying the comparison principle for nonlinear elliptic equations (cf. [GT, Theorem 17.1]) gives uniqueness.\qed

\it Counterexamples. \normalfont We give some examples to show that the inequality in Theorem 3.1 fails, if we do not assume any of the conditions on $\N.$ For small $\epsilon$, let $\Omega \subset \R^2$ be the annulus $$ \{x=(x_1, x_2)\, \mid \, |x|^2 \in [1,1+\epsilon]\}$$ and  $$\Gamma = \{(x_1, x_2, x_1, x_2) | x \in \Omega \} \subset \R^4$$ the gradient graph of $|x|^2/2.$  Let $\N$ = $\N_1 \cup \N_2$, where $$\N_1 =  \{(x_1, x_2, x_1, x_2) \,|\, x \in \overline{B_1}\} $$ and $$ \N_2 = \{(x_1, x_2, x_1+\eta(x), x_2+\eta(x))\, |\, x \in \overline {B_{1+\epsilon}} \} $$ where $\eta(x)$ is a small function which is positive on the interior of $B_{1+\epsilon},$ and vanishes on $\partial B_{1+\epsilon}.$   With a suitable orientation, the disconnected set $\N$ is homologous to $\Gamma$. We see that $ Vol(\N)$ is very close to $2\pi,$ whereas $ Vol(\Gamma) $ is very close to $0$, as we have chosen $\epsilon$ and $\eta$ small.  

To obtain a connected $\Sigma$ for which Theorem 3.1 fails, we alter the previous example slightly. With $\Omega$ as above, define a small ``bridge" region  $$\Omega' = \Omega \cap \{x_1 > 0\} \cap \{x_2 \in [-\epsilon, \epsilon] \} . $$  Let $\Gamma$ be the gradient graph of $|x|^2/2$ over $\Omega\remove\Omega',$ let $\N_1$ be the gradient graph of $|x|^2/2$ over $B_1 \cup \Omega',$ and let $\N_2$ be as above.  Then $\N = \N_1  \cup \N_2$ is connected and homologous to $\Gamma$, and is space-like except on a singular set which is one-dimensional, but nonetheless disconnects $\N$. Again we have $Vol(\N) $ close to  $2\pi$ and $Vol(\Gamma)$ close to $ 0,$ so the inequality in Theorem 3.1 does not hold.      
\bigskip

One can also approach the theory of special Lagrangian calibrations for pseudo-Eulidean space on the level of homogeneous spaces, following the approach of Harvey and Lawson.  First, noting that $U(n) = O(2n) \cap Sp(2n) $, we consider the group $G = O_{g_0}( n,n ) \cap Sp(2n)$, where $$O_{g_0}(n,n) = \{A \in M(2n, \R) | (Au,Av)_{g_0} = (u,v)_{g_0}, \, \forall u, v \in \Rnn \}.$$ Some standard Lie algebra computations will show that $$  G = \{ \begin{pmatrix} A & 0 \\ 0 & B  \\ \end{pmatrix} | AB^T = I \in M(n,\R) \}.$$  Next, consider the space $L^+$ of all space-like Lagrangian planes. Given $\xi \in L^+$, choose an oriented basis, $\xi_1,...,\xi_n,$ with each $\xi_i \in \Rnn,$ which is orthonormal with respect to $g_0$, and let $(A,B)$  be the $n \times 2n $ matrix with rows $\xi_i$. Then associate to each $\xi \in L^+$ the element $$\begin{pmatrix} A & 0 \\ 0 & B  \\ \end{pmatrix} \in G,$$ which is unique up to a factor of $SO(n),$ where the $SO(n)$ action on $G$ is defined by $$ S \cdot \begin{pmatrix} A & 0 \\ 0 & B  \\ \end{pmatrix} = \begin{pmatrix} SA & 0 \\ 0 & SB  \\ \end{pmatrix}  .$$  This gives the set of space-like Lagrangian planes $L^+$ the structure of a homogeneous space $L^+ \approx G / SO(n) .$ The homomorphism  $\Theta^+ : G \rightarrow O_{g_0}(1,1)$ given by $$ \begin{pmatrix} A & 0 \\ 0 & B  \\ \end{pmatrix} \mapsto \begin{pmatrix} \det(A)&0 \\0 & \det(B)   \\ \end{pmatrix} $$descends to a map $ L^+ \rightarrow O_{g_0}(1,1) \approx H^1$ where $H^1$ is the pseudo-circle $  \{ (s,t) |\,t >0 ,\, st = 1\} .$  A special Lagrangian submanifold is one whose tangent planes lie in a single fiber of $\Theta^+.$  We can see then that the calibrations are $$\Phi_c(\xi) = {\frac{1}{2}} ( c \det A + \frac{1}{c} \det B )  \geq Vol(\xi) $$ and the special Lagrangian equations are $\det D^2 u = c^2 > 0 ,$ where the value $c^2$ is the pseudo-phase analogous to the phase $c$ in section 2.

\bigskip
\section[]{A family of nonlinear equations}

For $ \thetat \in [0,\frac{\pi}{2}] $,  let $M_\thetat = (\R^n \times \R^n, g_\thetat),$ with $g_\thetat$ as defined in section 1.    For $\thetat < \frac{\pi}{4}, $ $M_\thetat$ is a pseudo-Euclidean space of index $n$. For $\thetat > \frac{\pi}{4},$ $M_\thetat$ is a Euclidean space.  For $\thetat = \frac{\pi}{4}, $ $M_\thetat$ carries a degenerate metric of rank $n.$  In any case, an extremal (minimal or maximal) submanifold $\N$ satisfies $ H= 0 $ along $\N$, where $H$ is the mean curvature vector along $\N$. 

 \begin{lemma}Suppose $u \in C^3(\Omega)$ and $ \Gamma = (x, \D u(x) ) $ defines an extremal surface in $M_\thetat,$ $\thetat \neq \pi / 4$ . If $D^2 u$ is diagonalized at a point $p=(x_0, \nabla(x_0)),$ the mean curvature satisfies
$$ H = 0 \overset{p} =  \sum_{i, k}  (\frac{u_{kii}}{\sin \thetat (1 + \lambda_i^2) + 2 \cos \thetat \lambda_i} \partial_{y_k})^N . $$ where $^N$ denotes orthogonal projection onto  $N(\Gamma)$, the normal bundle of $\Gamma$. 
\end{lemma}

\it Proof. \normalfont  For an embedded submanifold of $\Rnn$ given by $$ f:\Omega \rightarrow \Gamma \hookrightarrow M_\thetat, $$ we can compute the mean curvature vector $H$ by $$H = (g^{ij}f_{ij})^N $$  where  $g$ is the induced metric on $\Gamma,$ and $g^{ij} = (g^{-1})_{ij}.$ In this case, $$f(x) = (x_1, ... , x_n, u_1(x), ... u_n(x) ) ,$$ so
$$ H =  (g^{ij}(0,0,...0, u_{1ij}, ... u_{nij}))^N. $$ When $D^2 u = $ diag$(\lambda_1, ..., \lambda_n )$,  $$H =  \sum_{i, k}  (\frac{u_{kii}}{\sin \thetat (1 + \lambda_i^2) + 2 \cos \thetat \lambda_i} \partial_{y_k})^N = 0, $$ and $(\partial_{y_k})^N$ is a linear combination of $\partial_{x_k}$ and $\partial_{y_k}.$ For $ k = 1,...,n $ the vectors $(\partial_{y_k})^N $ are independent and form a basis for the normal space at $p.$  It follows that $$ \sum_{i}  \frac{u_{kii}}{\sin \thetat (1 + \lambda_i^2) + 2 \cos \thetat \lambda_i} = 0 $$ for all $k.$ $ \qed$
\bigskip

We now define the nonlinear operators $F^\thetat (D^2 u)$ by the equations (1.3) - (1.7).

At a point $p = (x_0, \nabla u(x_0))$, we may diagonalize $D^2 u $ and differentiate $F^\thetat(D^2 u ).$ If the eigenvalues $\lambda_i$ are all distinct, then each is differentiable, and we have 
$$ \frac{\partial}{\partial x_k} \lambda_i \overset{p}= \frac{\partial}{\partial x_k}u_{ii} = u_{iik} $$ and $$\frac{\partial}{\partial{x_k}} F^\thetat(D^2u(x))  = \ds\sum_i \frac{1}{\sin \thetat(1+\lambda_i^2) + 2 \cos \thetat \lambda_i} \frac{\partial}{\partial{x_k}} \lambda_i $$ \begin{equation} \overset{p} = \sum_i \frac{u_{kii}}{\sin \thetat (1 + \lambda_i^2) + 2 \cos \thetat \lambda_i} .\end{equation} If the eigenvalues are not distinct, then one can use an eigenspace projection argument to verify that $$ \frac{\partial}{\partial x_k} \sum_{i \in I_j} \lambda_i \overset{p}= \sum_{i \in I_j} u_{iik} $$ where $ I_j= \{i \in 1,...,n | \lambda_i \overset{p} = \lambda_j \}, $ so that (4.1) holds in this case as well.  It follows that the mean curvature $H$ vanishes if and only if $\nabla F^\thetat(D^2u) = 0$, that is, if and only if $F^\thetat(D^2 u )$ is constant. 

\bigskip

For this family of nonlinear equations,  $$\frac{\partial F^\thetat}{\partial \lambda_i} = g^{ii} $$  so we see that solutions are elliptic precisely when the Lagrangian submanifold is space-like. If we differentiate a second time with respect to $\lambda_i$,  we see that the equation is convex when all $\lambda_i < -a$,  and concave when all $\lambda_i > - a$.  We note that as $\thetat \rightarrow  \frac{\pi}{4}, $ equation (1.4) with 
$ c = n \ln (a-b) $ becomes (1.5) with $ c= \frac{n}{2}$.  Equation (1.5) has been studied by Flanders [F], who obtained Theorem 1.3 ii).

\section[]{Calibrations for $M_\thetat$ }

\bigskip

Gradient graphs for solutions to (1.4) and (1.6) give rise to calibrated submanifolds of $M_\thetat$, that are in fact isometric to calibrated submanifolds of $M_0$ and $M_\frac{\pi}{2}$. The calibrations for $M_\thetat$ may be obtained by pulling back the calibrations on $M_0$ and $M_\frac{\pi}{2}$  described in sections 3 and  2 via an isometry.

 Throughout this section we will again be using the constants $a = \cot \thetat $ and $b = \sqrt{|\cot^2 \thetat - 1| } ,$ for $\thetat  \in [0, \frac{\pi}{2}]$, as well as the constants defined by
\newcommand{\A}{{\sigma}}
\newcommand{\B}{{\tau}}
$$\A = \frac{\sqrt{\cos\thetat +\sin\thetat}+\sqrt{|\cos\thetat -\sin\thetat|}}{2}$$ 
$$\B = \frac{\sqrt{\cos\thetat +\sin\thetat}-\sqrt{|\cos\thetat -\sin\thetat|}}{2}.$$  

\bigskip
We start with the pseudo-Euclidean metrics.  For $\thetat < \frac{\pi}{4}$, the map  $  \varphi_\thetat: M_\thetat \rightarrow M_0,  $ represented by the $2n \times 2n$ matrix $$  \varphi_\thetat = 
\begin{pmatrix} \A I &  \B I \\ \B I &  \A I  \\ \end{pmatrix} $$is an isometry (up to a constant factor).  The Lagrangian condition is preserved under $\varphi_\thetat$, so this isometry maps Lagrangian $n$-surfaces  to Lagrangian $n$-surfaces.  Pulling back the calibrations on $M_0$ described in section 3, $M_\thetat$ becomes a calibrated manifold.  The isometry gives an equivalence between the homogeneous space structure for the space-like Lagrangian planes of $M_\thetat $ and the homogeneous space structure for the space-like Lagrangian planes of $M_0$ presented in section 3.  

\bigskip

\it Proof of Theorem 1.2 i). \normalfont Suppose $u(x)$ is a solution to (1.4). Let $\Gamma = (x, \D u(x)) \subset M_\thetat $ be the graph of $ \D u$. Our goal is to show that $\Gamma $ is isometric to a special Lagrangian graph $\hat\Gamma \subset M_0.$   Take $\hat{\Gamma} = \varphi_\thetat(\Gamma)\subset M_{0} $. At a point $p = (x_0, \nabla(x_0))  \in \Gamma$, the tangent space of $\Gamma$  can be described by the span of the vectors 
$$ \partial_i = \partial_{ x_i} + \sum_j u_{ij}\partial_{ y_j}, \quad i= 1,...,n.  $$ Take  $D^2 u $ to be diagonalized at $p$, and push forward. The tangent space $T_{\varphi_\thetat(p)}\hat{\Gamma} = (\varphi_\thetat)_*(T_p\Gamma) $ is the span of
$$ \hat{\partial_i} = (\varphi_\thetat)_*\partial_i = (\A + \B \lambda_i) \partial_{x_i}   + (\B + \A \lambda_i) \partial_{y_i} .$$ The space-like condition on $\Gamma$ imposes restrictions on the values of $\lambda_i,$ particulary,  $\lambda_i \neq -\A / \B ,$ so we may multiply each ${\hat\partial_i}$ by $  1 / (\A + \B \lambda_i)$ and see that $T_{\varphi_\thetat(p)}\hat{\Gamma}$ is spanned by
$$ \hat\partial_{x_i} + \frac{\B + \A \lambda_i}{\A + \B \lambda_i} \hat\partial_{y_i}.$$ The image $\hat\Gamma$ is a Lagrangian submanifold of $M_0$, so arises locally as gradient graph $\hat \Gamma = (\hat x, \D \hat u(\hat x)).$ From the above expression, the eigenvalues of $D^2\hat u$ are given by  
\begin{equation} {\hat\lambda_i} = \frac{\B  +   \A \lambda_i}{ \A  + \B \lambda_i} =\frac{  ({\lambda_i} + \B/\A ) }  {({\lambda_i} + \A/\B)}\frac{(\A)}{(\B)} . \end{equation} The $\lambda_i$'s satisfy (1.4), so noting that $ a + b = 1/(a-b) = \A/\B,$ we have
$$\prod  \frac{  {\lambda_i} + \B/\A  }  {{\lambda_i} + \A/\B} =   e^c > 0\,,  $$ and may conclude that $\hat u(\hat x)$ satisfies $$\prod \hat\lambda_i =  [\frac{\A}{\B}]^n \,e^c > 0 ,  $$ that is, $\hat u(\hat x)$ satisfies the Monge-Amp\`ere equation (1.1). It follows that $\hat \Gamma$ is a calibrated submanifold. The property of being calibrated is local and is preserved under isometries of the ambient manifolds, so $\Gamma$ is calibrated.  Theorem 1.2 i) then follows by the same reasoning as in the proof of Theorem 3.1. \qed 
\bigskip


The function $\hat u$ obtained above is only local.  In order to study solutions to (1.4) further, in particular to obtain uniqueness for Theorem 1.2 i)  and the Bernstein-type result, we transform $u$ into a solution to (1.1) which describes $\hat\Gamma$ globally, when possible.  
 
The isometry $\varphi_\thetat$ acts via $$\varphi_\thetat(x,y) = (\sigma x + \tau y,\tau x +\sigma y) = (p(x,y), q(x,y)) $$ so the image  $\hat\Gamma = \varphi_\thetat(\Gamma)$ lies in the set $p(\Gamma) \times q(\Gamma) \subset \R^n \times \R^n.$  As $\Gamma$ is parameterized by $\Omega,$ the isometry $\Gamma \rightarrow \hat\Gamma$ locally amounts to a change of coordinates $\Omega \rightarrow p(\Gamma),$ wherever $Dp$ is invertible.  If $p|_\Gamma$ is a bijection, then $\hat \Gamma$ is globally parameterized by $p(\Gamma)$. In the case where $p|_\Gamma$ is a bijection, $\hat \Gamma$ is a graph over $p(\Gamma)$ of the function $$ r = q \circ p^{-1} : p(\Gamma)
 \rightarrow \R^n. $$ Using $\Omega$ as coordinates for $\Gamma$  
$$p(x) = \A x +  \B \nabla u (x) $$ $$ q(x) = \B x + \A \nabla u (x) $$ hence
$$ Dp = \A I + \B D^2 u $$ $$ Dq = \B I + \A D^2 u$$ so $p$ is locally invertible if and only if no eigenvalues of $D^2 u $ attain the value $- {\A}/{\B} = -(a+b), $ which is precluded by the space-like condition.   Diagonalizing $D^2 u$, we see that $Dr$ is symmetric, so $r = \nabla \hat u$ for some $\hat u,$ and the eigenvalues of $D^2 \hat u$ become
$$ \hat \lambda_i = \frac{\B + \A\lambda_i}{\A + \B\lambda_i}  $$  as in (5.1).  Now if  $D^2 u \geq -a > -(a + b) $, then $Dp >  0$ and $p$ is injective on $\Omega.$ The function $\hat u$ is then a solution to the Monge-Amp\`ere equation defined on all of $p(\Gamma)$.  Further, the inequality is uniform,  $Dp \geq \epsilon > 0. $ Thus if $\Omega = \R^n,$ $p$ is a bijection on $\R^n$, so $\hat u$ will then be a convex solution to (1.1) on all of $\R^n.$  Theorem 1.3 i) follows from the famous result

\begin{theorem}[Calabi \mbox{[C]}, Pogorelov \mbox{[P]}] Any convex solution to (1.1) on all of $\R^n$  is a quadratic polynomial. 
\end{theorem} 
Using a recent, more general result, we see that we may drop the restriction that $D^2 \geq - a .$

\begin{theorem}[Jost-Xin, \mbox{[JX, Theorem 4.2]}] Let M be a space-like extremal $m$-surface in $\R^{n+m}_n.$ If M is closed with respect to the Euclidean topology, then M must be a linear subspace. 
\end{theorem} 
\bigskip

\it Uniqueness. \normalfont From the above discussion, if we assume that $D^2u \geq -a,$ the surface $\Gamma$ can be described globally as the gradient graph of a solution to the Monge-Amp\`ere equation. In this case, that $\Gamma$ is the unique surface with this volume follows from Theorem 3.1.  The proof of Theorem 3.1 relied heavily on the fact that the space-like condition forced the projection $\Gamma \rightarrow \R^n_x$ to be an open map.  This is not the case for gradient graphs of solutions to equation (1.4), so the corresponding proof of uniqueness will not work.  We do not have a proof of uniqueness for equation (1.4) with no restrictions on $D^2 u ,$ nor for the analogous result for solutions to (1.6).  As in [HL, Theorem 5.8], we do have uniqueness for solutions to (1.4) and (1.6) whenever the boundary data is analytic, as an application of the Cauchy-Kowaleswki Theorem. 

\bigskip

\it Proof of Theorem 1.2 ii) \normalfont  Let $ P: \Rnn \rightarrow \R^n$ be the map 
$$ (x,y) \mapsto \frac{x+y}{2}. $$  Given the degenerate metric $g_\frac{\pi}{4}$ on $\Rnn$, this map is an ``isometry" in the sense that $g_\frac{\pi}{4}= P^*\delta_0$ . The graph $\Gamma = (x, F(x))$ is isometric to $P(\Gamma)$ for any $F(x)$ such that $DF $ avoids the eigenvalue $\lambda_i = -1$. Let $\N$ be any space-like surface with boundary $\partial \N =  \Gamma \cap \{x \in \Omega\}. $  Then if $\omega = x_ndx_1\wedge...\wedge dx_{n-1} ,$ so that $dw = dVol$ on $\R^n,$ $$Vol(\N) = \int_\N dVol = \int_\N P^* dVol = \int_\N P^*dw = $$ $$\int_\N dP^* w = \int_{\partial \N} P^*w =\int_{\partial\Gamma}  P^*w = Vol(\Gamma). \,\,\,\qed$$

\bigskip

\bigskip 

\it Proof of Theorem 1.2 iii). \normalfont 
As in the pseudo-Euclidean case, for  $\thetat > \frac{\pi}{4}$, the map  $  \varphi_\thetat: M_\thetat \rightarrow M_\frac{\pi}{2}, $ is an isometry.  For a given $n$-plane $\xi$, 
$$Vol_{g_\thetat}^2(\xi) =Vol_{\delta_0}^2(\varphi \xi) \geq Vol(\varphi \xi \wedge J\varphi \xi ) $$
$$= \det_\R A = ||\det_\C A||^2 = \alpha^2(\xi) + \beta^2(\xi) \geq \alpha^2(\xi) ,$$ where $A \in M(n,\C) \subset M(2n, \R) $ is the map sending $\partial_{x_i} \mapsto \varphi \xi_i$, extended by complex linearity as before.  By Hadamard's Inequality, we have equality if and only if $\xi$ is Lagrangian and $\beta(\xi)=0.$

If $\xi$ is Lagrangian, then $\varphi \xi_1,...,\varphi \xi_n$ is an orthonormal basis for $\C^n$, so $A \in U(n).$ Using the map
$$\Theta: L \approx U(n) / SO(n) \overset{\det_\C}\rightarrow S^1.$$ the special Lagrangian equations (1.6) can be deduced from
$$ \arg \det  (\A + \B D^2u + \I (\B +\A D^2 u)) = \sum \arctan (  \frac{\B  +   \A \lambda_i}{ \A  + \B \lambda_i}) =   c.  $$ If follows that any gradient graph $\Gamma = (x, \nabla u(x))$ for $u(x)$ satisfying this equation is calibrated, and therefore an absolutely volume minimizing submanifold.  \qed

\bigskip

By differentiating with respect to $\lambda$, one can verify the identity
$$\sum  \arctan \frac{\lambda_i + a }{b}  + C_\thetat  = \sum \arctan (  \frac{\B  +   \A \lambda_i}{ \A  + \B \lambda_i} ) $$ where $$ C_\thetat =  \arctan ( \frac{ \B}{\A})  - \arctan{ \frac{a}{b} }.$$ It follows that a solution to (1.6) also satisfies 
$$  \sum  \arctan \frac{\lambda_i + a }{b}  = c - C_\thetat,  $$ and the function $$ v(x) = \frac{u(x)}{b}+ \frac{a}{2b}|x|^2 $$ is a solution to (1.2).  If $v(x)$  is convex, which is the case if $ D^2 u \geq -a $, we may conclude that $v(x)$  is a quadratic polynomial by applying the following result.
\begin{theorem}[Yuan \mbox{[Y1]}] Suppose $u \in C^2(\R^n)$ is a convex solution to (1.2).  Then $u(x)$ is a quadratic polynomial.  \end{theorem}

Similarly, any solution to (1.6) with $c > {(n-2)\pi}/ 2 + C_\thetat$ is a quadratic polynomial, by the theorem

\begin{theorem}[Yuan\mbox{ [Y2]}] Suppose $u \in C^2(\R^n)$  is a solution to (1.2) with $c>\frac{(n-2) \pi}{2}.$ Then $u(x)$ is a quadratic polynomial.  

\end{theorem}

\bigskip

\section[]{Example}
The equations (1.3)-(1.7) can be manipulated to take the form $$f(\lambda_1,...,\lambda_n) = 0,$$ where $f$ is a polynomial of degree no higher than $n.$  For the equation (1.4) with $c= 0$, the polynomial $f$ has degree $n-1.$ When $n=2$, (1.4) with $c=0$ becomes a linear equation, taking the form $\Delta u = -2a.$   We exploit this degeneracy in order to write down an explicit solution to (1.4).

Fix some $\thetat \in (0,\frac{\pi}{4}) $, and define $u(x_1,x_2)$ on $\Omega = \{ x_1>0 \subset \R^2 \}$

$$ u(x_1,x_2) = -\frac{a}{2}(x_1^2+x_2^2) + ke^{x_1} \cos x_2, $$for $k$ some large constant, and $a = \cot \thetat$.  The eigenvalues of $D^2 u $ are
$$\lambda_1 = -a - ke^{x_1} , $$
$$\lambda_2 = -a + ke^{x_1} , $$ Clearly, $u$ is a solution to (1.4) with $ c = 0$, that is, $\Delta u = -2a,$ and one can check that the resulting surface $\Gamma = (x, \nabla u(x))$ is space-like in $(\R^2 \times \R^2, g_\thetat)$.  We transform this to a maximal submanifold of $(\Rnn, g_0),$ using $\varphi_\thetat.$  This isometry  acts via $$\varphi_\thetat(x,y) = (\sigma x + \tau y,\tau x +\sigma y)$$ so the image $ \varphi_\thetat(\Gamma) = \hat\Gamma \subset (\R^4, g_0) $ lies in $p(\Gamma) \times \R^2$, where $ p(x,y) = \A x + \B y. $ Parameterizing $\Gamma$ by $\Omega$, the map $p|_\Gamma$ takes the form $$ p(x_1,x_2) = \A \cdot (x_1,x_2) + \B \cdot \nabla u(x_1,x_2),$$ $$ \mbox{or}\,\, p(z) = (\A-\B a )z + \B ke^{\bar{z}}, $$ where $z = x_1 + ix_2.$  For very large $k$, a simple volume computation over a vertical strip shows that $p$ can not be injective on $\Gamma$, so $\hat\Gamma$ will not be described globally as a graph over $p(\Gamma)$.   This surface $\hat\Gamma$ is a maximal Lagrangian surface which is globally described as a solution to (1.4) and locally described as a solution to (1.1), but which is not globally described as a solution to (1.1).  This surface is not complete, as we have restricted the domain to a half-plane, but by the result of Jost and Xin [J-X, Theorem 4.1] one would not expect to find a complete surface with this property.


\begin{thebibliography}{9999}
%


\bibitem[C] {C}E. Calabi. \it Improper Affine Hyperspheres of convex type and a generalization of a theorem by K. J\"orgens. \normalfont Michigan Math J. \bf 5\normalfont(2)   (1958), 105-126. 


\bibitem[F] {F}H. Flanders, \it On certain functions with positive definite Hessian, \normalfont  Ann. of Math. (2)  \bf 71 \normalfont 1960 153--156.


\bibitem[GT] {GT}D. Gilbarg and N. Trudinger, \it Elliptic partial differential equations of second order.    \normalfont Springer-Verlag, Berlin, 2001. 
\medskip

\bibitem[H] {H}N. J. Hitchin. \it The moduli space of 
special Lagrangian submanifolds.  Dedicated to 
Ennio De Giorgi. \normalfont Ann. Scuola Norm. 
Sup. Pisa Cl. Sci. (4) \bf 25\normalfont 
(1997), no. 3-4, 503--515 (1998). 




\medskip

\bibitem[HL] {HL}R. Harvey and H. B. Lawson. Jr., \it Calibrated Geometries, \normalfont Acta Math \bf 148 \normalfont (1982), 47-157. 

\medskip





\bibitem[JX] {JX}J.Jost and Y.L.Xin, \it  Some aspects of the global geometry of entire space-like submanifolds,  \normalfont
 Result. Math. \bf 40 \normalfont(2001) 233--245. 

\medskip



\bibitem[L] {L}J. Lee, \it Introduction to smooth manifolds, \normalfont Springer-Verlag, New York, 2003. 
\medskip



\bibitem[M] {M}R.C. McLean, \it 
Deformations of calibrated submanifolds, 
\normalfont

Comm. Anal. Geom. \bf 6 \normalfont (1998), 
no. 4, 705--747.

\medskip


\bibitem [O] {O} R. Osserman \it Minimal Varieties \normalfont   Bull. Amer. Math. Soc.  \bf 75 \normalfont (1969),  1092-1120.


\medskip

\bibitem [P] {P}A.V. Pogorelov,\it On the improper convex affine hyperspheres , \normalfont   Geometriae Dedicata \bf 1 \normalfont (1972),  33--46.

\medskip


\bibitem[Y1]{Y1} Y. Yuan, \it A Bernstein problem for special Lagrangian equations \normalfont
Invent. Math. \bf 150\normalfont (2002), no. 1, 117--125.

\medskip

\bibitem[Y2]{Y2} Y. Yuan, \it Global Solutions to special Lagrangian Equations, \normalfont Proc. Amer. Math. Soc.\bf  134  \normalfont (2006),  no. 5, 1355--1358. 

\medskip


 


\end{thebibliography}
\end{document}